 \documentclass[12pt,english]{article}

\usepackage{graphicx,psfrag,epsfig,amsfonts,amssymb,amsmath,babel}

\usepackage{color}

\newtheorem{definition}{Definition}[section]

\newtheorem{lemma}[definition]{Lemma}
\newtheorem{theorem}{Theorem}

\begin{document}

\title{Oscillating statistics   of transitive dynamics}

\author{Eleonora Catsigeras   \\ Instituto de Matem\'{a}tica y Estad\'{\i}stica \lq\lq Prof. Ing. Rafael Laguardia\rq\rq, \\
 Facultad de   Ingenier\'{\i}a, \\ Universidad de la Rep\'{u}blica, \\  URUGUAY \\
 E-mail: eleonora@fing.edu.uy}

\date{June 12, 2015}

\maketitle

\begin{abstract}
   We   prove that   topologically generic orbits   of   $C^0$, transitive  and non-uniquely ergodic dynamical systems, exhibit an extremely oscillating  asymptotical statistics. Precisely, the minimum weak$^*$ compact set of invariant probabilities that describes the asymptotical statistics of each   orbit of a residual set  contains all the ergodic probabilities. If besides $f$ is ergodic with respect to the Lebesgue measure, then also Lebesgue-almost all the orbits exhibit that kind of extremely oscillating statistics.

\end{abstract}


\section{Introduction}
\label{SectionIntroduction} \vspace{-4pt}

We will study the statistical average for typical orbits of transitive dynamics, under a non traditional viewpoint.

On the one hand, the traditional viewpoint studies the limit in the future of the Birkhoff averages, \em   starting always from the same initial  point, \em and for Lebesgue-positive sets of orbits in the future.
So, under this traditional viewpoint,  the \lq\lq statistics\rq\rq \ of the system (at least for $C^2$-dynamical systems  with some kind of hyperbolicity), is mainly obtained
from the existence of physical measures,
  of Sinai-Ruelle-Bowen (SRB) measures, and of     Gibbs measures (see for instance the survey  \cite{Buzzi}).

 Relevant advances
on the study of the asymptotic behavior of the time-averages   from the traditional viewpoint
can be found  for instance
in the following articles. In  \cite{Viana-Yang} Viana and Yang study the existence of
physical measures   for partially
hyperbolic systems with one-dimensional center
direction.   Bonatti's survey \cite{Bonatti}   gives
an overview of the state of art in the theme
of
the  asymptotical dynamics of $C^1$-differentiable systems
     from the topological viewpoint.
In \cite{Liverani} Liverani proves
   that piecewise $C^2$ expanding
maps may exhibit Gibbs measures without needing the bounded distortion
property.

On the other hand, instead of adopting the traditional viewpoint, along this paper we will study the-time averages that   start at \em any future iterate of the initial point. \em This viewpoint is based on a philosophical argument: the way that   the observers  in the future will perceive the   forward statistics of the  system, is not   the way that it is computed  today. In fact, today the observers   compute the Birkhoff average  along the finite future piece   orbit of length $n$ (which we like to call \lq\lq the clima\rq\rq), by the mean value of the observable functions from time 0 to $n$. But   the observers in the future -who will live, say, at time $m >0$- will compute their Birkhoff average along the finite piece  of orbit of  length $n$ (i.e. they will perceive their clima),    by the mean value of the observable functions between time  $m$ and time $m + n$.

This non-traditional viewpoint of studying the Birkhoff averages and their limits (i.e. the statistics)  does not  give  preferences to different initial observation instants. So, our conclusions include also the prediction of \lq\lq all the climas \rq\rq \ that   the observers in   the future  will perceive.

The key result is Theorem \ref{theorem0}:

 \em Topologically typically,  the clima observed at infinitely many times in the future must widely differ from the clima observed at present time, provided that the dynamics is   deterministic   (non hazardous), transitive and non-uniquely ergodic. \em

This is   an unexpected result, taking into account that the system is \em autonomous and deterministic. \em
Nevertheless, the idea of the proof of Theorem \ref{theorem0} is extremely simple. The route of its proof is the result of joining the following three simple observations.
First, if the system is transitive, then its topologically generic orbits in the future are dense. Second, for any ergodic measure $\mu$, and for any $\mu$-typical point $x_0$, the Birkhoff average starting at $x_0$ converges to $\mu$. So, for any $\epsilon >0$, for any \em fixed \em $n$ sufficiently large, and for any point $x$ close enough $x_0$,  the Birkhoff average starting at $x$  is  $\epsilon$-near $\mu$.Third, any dense orbit   in the future has such an  iterate $x$ close enough $x_0$.

Thus, one concludes that the Birkhoff averages, with fixed $n$ but starting at different  points in the future of the same orbit, oscillate among \em all \em the ergodic measures of $f$, when $n \rightarrow + \infty$.

Even if the main theorem is the consequence of the latter simple observations, and no more proof than the above argument would be needed, we will include  all the details of this proof   (see Section \ref{SectionProofs})  to be readable by a wide class of scientists and students.

\subsection {Mathematical background}

Let $M$ be a compact manifold of finite dimension. Let $f: M \mapsto M$ be   continuous. We consider the dynamical system obtained by   iteration of $f$ in the future, i.e. the family of orbits $\{f^n(x)\}_{n \in \mathbb{N}}$ with initial condition  $x \in M$. This  dynamical system is composed by the  solutions of the recurrent equation $x_{n+1}   = f(x_n)$.

We denote by ${\mathcal M}$ the space of all the probability measures in $M$, endowed with the weak$^*$ topology (see for instance Definition 6.1 of \cite{libroteoergodica}). That is, if $\mu_n$ is a sequence of probability measures in $M$, we define \begin{equation}
\label{eqnWeak*Topology} \lim \mu_n = \mu \in {\mathcal M} \mbox{ if and only if }
\lim _{n \rightarrow + \infty} \int \varphi \, d \mu_n = \int \varphi \, d \mu \ \ \ \forall \ \varphi \in C^0(M, \mathbb{R}),\end{equation}
where $C^0(M, \mathbb{R})$ is the space of continuous real functions in $M$, with the supremum norm.

Recall that a  measure $\mu \in {\mathcal M}$ is invariant by $f$ if $\mu(f^{-1}(A)) = \mu(A) $ for any Borel-measurable set $A \subset M$. We denote by ${\mathcal M}_f \subset {\mathcal M}$ the space of $f$-invariant probability measures, and by ${\mathcal E}_f \subset {\mathcal M}_f$ the set of ergodic probability measures for $f$.

We recall that $\mu \in {\mathcal M}_f $ if and only if
$$\int \varphi \, d \mu = \int \varphi \circ f \, d \mu \ \ \ \forall \ \varphi \in C^0(M, \mathbb{R}).$$ (See for instance  Theorem 6.8 of \cite{libroteoergodica}).

To each initial state $x \in M$, or equivalently  to each orbit $\{f^n(x)\}_{n \in \mathbb{N}}$,  we associate the double-indexed  sequence $\{\sigma_(m, n)(x)\}_{(m, n) \in \mathbb{N}^2}$ of non necessarily invariant  probability measures $\sigma_{m,n}(x)$, which we call \em empirical probabilities\em,  defined by:
\begin{equation}\label{eqn00}\sigma_{ m,n }(x) = \frac{1}{n} \sum_{j= m}^{m + n-1} \delta_{f^{j } (x)},\end{equation}
where $\delta_y$ is the Dirac-Delta probability measure supported on the point $y \in M$. In other words, the empirical probability $\sigma_{m, n}(x)$ is the probability distribution that is observed during  a statistical experiment  on which one computes the Birkhoff average (i.e. the temporal average) of the observable  functions  $\varphi:M \mapsto \mathbb{R}$ along a finite piece of the orbit of $x$, from time $m$ to time $m + n-1$. Precisely:
\begin{equation}
\label{eqn01}
\frac{1}{n} \sum_{j= m}^{m + n-1} \varphi (f^{j }(x))  = \int \varphi \, d \sigma_{m, n} (x).\end{equation}
We agree to call the double-indexed sequence
$\{\sigma_{m,n}(x)\}_{(m,n) \in \mathbb{N}^2}$ of empirical probabilities \em  the complete   future  statistics \em of the orbit of $x$. For the sake of concision we call it the \em statistics \em of $x$.

Since the space ${\mathcal M}$ is metrizable and weak$^*$-compact,  it is sequentially compact (see for instance Theorems 6.4 and 6.5 of \cite{libroteoergodica}). Thus, any   sequence $\{\sigma_{m(n), n}(x)\}_{n \in \mathbb{N}}$ of empirical probabilities  has convergent subsequences when $n \rightarrow + \infty$.

We agree to call the set of all limit probability measures  of all   such sequences of empirical probabilities in ${\mathcal M}$, \em the asymptotical statistics \em of the orbit of $x$ (Definition \ref{definitionAsymptoticalStatistics}).

%
%

\subsection{Statement of the results}

If $f:M \mapsto M$ preserves the Lebesgue measure $m$ and is ergodic, then the sequence $\{\sigma_{0,n}(x)\}_{n \in \mathbb{N}}$   is convergent  for  Lebesgue-almost all $x \in M$ (see for instance Theorem 6.12 (ii) of \cite{libroteoergodica}). In other words, its limit set is a singleton. Also, if there exists a unique physical measure whose basin of statistical attraction covers Lebesgue almost all the points, or if there exists a unique SRB-like measure, then the limit set of the sequence $\{\sigma_{0,n}(x)\}_{n \in \mathbb{N}}$ is   a singleton for Lebesgue-almost all $x \in M$ (see   \cite{polaca}).

In contrast, if  instead of restricting to the case $m(n)= 0$, we consider all the sequences of the form $\{\sigma_{m(n),n}(x)\}_{n \in \mathbb{N}}$ where $m: \mathbb{N} \mapsto \mathbb{N}$, then the limit set may be non convergent, and moreover, extremely oscillating (see Definition \ref{definitionExtemeOscilAE}). In fact, in this paper we prove the following result:

\begin{theorem}
  \label{theorem1}

  Let $f: M \mapsto M$ be continuous, preserving the Lebesgue measure of $M$ and ergodic with respect to it, but non uniquely ergodic. Then   Lebesgue-almost all the orbits of $f$ have  extremely oscillating asymptotical statistics. Precisely, it contains   all the ergodic probability measures of $f$.
  \end{theorem}

Let us state a similar result that holds for maps that do not preserve the Lebesgue measure.
In Theorem 3.6 of \cite{AbdenurAndersson}, Abdenur and Andersson  studied    the limit set of the sequence $\{\sigma_{0,n}(x)\}_{n \in \mathbb{N}}$  for  Lebesgue-almost all the orbits of $C^0$-generic  maps. Such generic systems do not preserve the Lebesgue measure.   They proved that  the particular sequence  $\{\sigma_{0,n}(x)\}_{n \in \mathbb{N}}$ of empirical probabilities  is convergent for Lebesgue-almost all $x \in M$. So, its limit set is a singleton.

 Now,   for transitive and non-uniquely systems, we observe all the sequences $$\{\sigma_{m(n),n}(x)\}_{n \in \mathbb{N}} \mbox{ with any } m: \mathbb{N} \mapsto \mathbb{N},$$ instead of restricting to the case $m(n)= 0$. Let us apply  a topological criterium   instead of  a Lebesgue-probabilistic criterium when selecting the \lq\lq relevant\rq\rq \ orbits of the system. With such an agreement, we say that an orbit is generic if it belongs to a residual set in $M$. Then the asymptotical statistics is far from being a singleton: it is extremely oscillating. In fact,   we prove the following  result:

  \begin{theorem}
  \label{theorem0}

  Let $f: M \mapsto M$ be continuous, transitive and non uniquely ergodic. Then generic orbits of $f$ have  extremely oscillating asymptotical statistics. Precisely, any ergodic probability for $f$ belongs to the asymptotical statistics of each generic orbit.
  \end{theorem}

Theorems \ref{theorem1} and \ref{theorem0} imply the necessary extremely changeable   \lq\lq clima\rq\rq,   i.e. the time averages of the observable functions along   finite pieces of all the relevant  orbits in the ambient manifold $M$ vary so much  in the long term, to  approach all the extremal invariant   probabilities of the system (the ergodic measures).    Even if the system is fully deterministic and it is governed by an autonomous and unchangeable recurrence equation, even if the parameters in this equation are fixed, even if the states along the deterministic orbit are not perturbed, no topologically relevant orbit of the system has a predictable statistics along its long-term future evolution. On the contrary, its asymptotical statistics is extremely changeable in the long-term  future, exhibiting at least, as many probability distributions as ergodic measures of $f$ exist.

\vspace{.3cm}

 This paper is organized as follows: In Section \ref{SectionDefinitions&Statements} we state the precise mathematical definitions to which the results refer, and in Section \ref{SectionProofs} we include the proofs of Theorems \ref{theorem1} and \ref{theorem0}.

\section{Definitions}
\label{SectionDefinitions&Statements} \vspace{-4pt}

Since the double-indexed sequence of empirical probabilities $\{\sigma_{m,n}(x)\}_{(m,n) \in \mathbb{N}^2}$ completely describes de statistics (i.e. the time-average) of any finite piece of the orbit of $x$, the  limit set $as  (x)$ in the space of probabilities describes what we call the   asymptotical statistics   of the orbit, according to the following definition:

\begin{definition}
\label{definitionAsymptoticalStatistics} \em {\bf (Asymptotical statistics $as(x)$ in the space of probabilities)}

The \em asymptotical statistics  of the orbit of $x \in M$, \em which we denote by $as (x)$,  is the set composed by all the limits in ${\mathcal M}$ of the convergent subsequences of any sequence $\{\sigma_{m(n), n} (x)\}_{n \in \mathbb{N}}$ of empirical probabilities of $x$, where $m: \mathbb{N} \mapsto \mathbb{N} $ is any mapping from the set of natural numbers to itself. Precisely:
\begin{equation} \label{eqn10}  as(x):= \{\mu \in {\mathcal M}: \ \exists \ \{m_i, n_i\}_{i \in \mathbb{N}} \subset \mathbb{N}^2 \mbox{ such that } n_i \rightarrow + \infty \mbox{ and }   \lim_{i \rightarrow + \infty} \sigma_{m_i, n_i}(x) = \mu\} \end{equation}
Following   the classical Krylov-Bogolioubov construction of invariant probabilities (see for instance the proofs of Theorems 6.9, \ 6.10, \ and Corollary 6.9.1 of \cite{libroteoergodica}), it is standard to check that:
$$as(x) \neq \emptyset, \ \ as(x) \mbox{ is weak$^*$-compact, and } as (x) \in {\mathcal M}_f \ \ \forall \ x \in M.$$ In other words, the asymptotical statistics  of $x$   is a nonempty compact set of probability measures which are invariant by $f$.

\end{definition}

  \begin{definition} \label{definitionConvergentAS} {\bf (Convergent or oscillating asymptotical statistics)} \em

  The orbit $\{f^n(x)\}_{n \in \mathbb{N}}$ is \em statistically convergent \em if its asymptotical statistics is composed by a unique probability measure, i.e.
  $$\# (as(x)) = 1.$$
  It is \em statistically oscillating \em if it is non convergent.

  \end{definition}

   We recall that $f$ is called \em uniquely ergodic \em if $\#{\mathcal E}_f = 1$ (see for instance \cite{Mane}).

\begin{definition} \label{definitionExtemeOscilAE} {\bf (Extremely oscillating asymptotical statistics)} \em

When $f$ is non-uniquely ergodic we say that the orbit $\{f^n(x)\}_{n \in \mathbb{N}}$ is \em statistically extremely oscillating \em if  its asymptotical statistics contains all the $f$-invariant ergodic probability measures. Namely:
  $$ as(x)  \supset{\mathcal E}_f , \ \  \#{\mathcal E}_f > 1.$$

   \end{definition}

  \begin{definition}
  \label{definitionTransitivo}{\bf (Transitive system)} \em The dynamical system by iterates of $f: M \mapsto M$ is called \em transitive \em if for any pair $(U,V)$ of nonempty open sets in $M$ there exists a positive iterate of $U$ that intersects $V$.

  Let us denote ${\mathcal T}(M)$ to the topology of $M$, i.e. the family of all the open sets of $M$. So, $f:M \mapsto M$ is transitive, by definition, if
  $$\forall \ (U,V) \in  {\mathcal T}(M) \times {\mathcal T}(M) \ \mbox{ if } U,V  \neq \emptyset, \mbox{ then } \ \exists \  n \in \mathbb{N}^+ \mbox{ such that } f^n(U) \cap V \neq \emptyset, $$
  where $ \mathbb{N}^+$ denotes the set of positive integer numbers.  Equivalently, $$\forall \ (U,V) \in  {\mathcal T}(M) \times {\mathcal T}(M) \ \mbox{ if } U,V  \neq \emptyset, \mbox{ then } \ \exists \  n \in \mathbb{N}^+ \mbox{ such that }   f^{-n}(V) \cap U \neq \emptyset.$$
  Recall that $M$ is a finite dimensional manifold. So,  $f:M \mapsto M$ is transitive if and only if there exists $x \in M$ whose   orbit in the future is dense in $M$.
  \end{definition}

  \begin{definition}
  \label{definitionGenericOrbit} {\bf (Residual sets and generic orbits)} \em

  According to Baire-category theory  a  set ${  R} \subset M$ is said \em residual  \em if it contains a countable intersection of open and dense subsets of $M$.  It is standard to check that the countable intersection of residual sets is residual. Since $M$ is a compact manifold, any residual set $R$ is dense, but not all dense sets are residual.

  Given a residual set $R \subset M$ we say that the orbits $\Big\{\{f^n(x)\}_{n \in \mathbb{N}}: x \in R\Big\} $ are \em generic. \em
  \end{definition}

  \section{The proofs} \label{SectionProofs}

   The weak$^*$ topology of the space ${\mathcal M}$ of probability measures is metrizable   (see for instance Theorem 6.4 of \cite{libroteoergodica}). We choose and fix a weak$^*$-metric in ${\mathcal M}$, which we denote by $\mbox{dist}$.

   \vspace{.3cm}

    To prove Theorems \ref{theorem1} and \ref{theorem0} we first state the following lemmas:

    \begin{lemma}
    \label{lemma0a} Let $f: M \mapsto M$ be continuous and ${\mathcal M}_f$ denote the space of $f$-invariant probability measures. Let $\mu \in {\mathcal M}_f$. Let $x \in M$ and let $as(x) \subset {\mathcal M}_f$ be the asymptotical statistics of the orbit of $x$, according to Definition \em \ref{definitionAsymptoticalStatistics}. \em Then $$\mu \in as (x) \ \mbox{ \em  if and only if  } \  x \in A(\mu), \mbox{ \em where } $$
    \begin{equation}
    \label{eqn11}
   A(\mu)= \bigcap_{\epsilon >0} \bigcap_{N \geq 1} \bigcup_{m \geq 0} \bigcup_{n \geq N} f^{-m}\Big(\{y \in M: \mbox{\em dist} (\sigma_{0,n}(y), \mu) < \epsilon\}    \Big) \end{equation}
    \end{lemma}

    \noindent{\em Proof:} From equality (\ref{eqn10}), $\mu \in as(x)$ if and only if:
    $$  \lim_{i \rightarrow + \infty} \sigma_{m_i, n_i}(x) = \mu $$
    for some sequence  $ \{m_i, n_i\}_{i \in \mathbb{N}} \subset \mathbb{N}^2 $   such that $ n_i \rightarrow + \infty $. This condition holds if and only if  for any $\epsilon >0$ there exists $i_0 \in \mathbb{N}$ such that
    $$\mbox{dist}(\sigma_{m_i, n_i}(x), \mu) < \epsilon \ \ \forall \ i > i_0.$$
    Since $\lim n_i = + \infty$, for any $N \geq 1$ there exists $i_1 \in \mathbb{N}$ such that $n_i \geq N$ for all $i > i_1$. We deduce that
     $$ n_i \geq N \mbox{ and } \mbox{dist}(\sigma_{m_i, n_i}(x), \mu) < \epsilon \ \ \forall \ i > \max\{i_0, i_1\}.$$
     In other words, $\mu \in ae(x)$ if and only if for all $\epsilon >0$ and all $N \geq 1$, the point $x $ belongs to the set
     $$\bigcup_{m \geq 0}\bigcup_{n \geq N} \{x \in M: \ \mbox{dist} (\sigma_{m,n}(x), \mu) < \epsilon \}.$$
     From equality (\ref{eqn00}) note that
     $$\sigma_{m,n}(x) = \sigma_{0,n}(f^m(x)).$$
     Then
     $$\mbox{dist}(\sigma_{m,n}(x), \mu) < \epsilon \ \Leftrightarrow \ \mbox{dist}(\sigma_{0,n}(f^m(x)), \mu) < \epsilon \ \Leftrightarrow $$ $$ x \in f^{-m}\Big(\{y \in M: \  \mbox{dist}(\sigma_{0,n}(y), \mu) < \epsilon\}\Big). $$
     We have proved that, $\mu \in a e (x)$ if and only if for all $\epsilon >0$ and all $N \geq 1$, the point $x $ belongs to the set
     $$\bigcup_{m \geq 0}\bigcup_{n \geq N} f^{-m}\Big(\{y \in M: \ \mbox{dist} (\sigma_{0,n}(y), \mu) < \epsilon \}\Big).$$
     We conclude that $x \in A_{\mu}$, where the set $A_{\mu}$ is defined by equality (\ref{eqn11}), ending the proof. \hfill $\Box$

    \begin{lemma}
    \label{lemma0b} If $\mu$ is ergodic, then for all $\epsilon >0$ there exists $N \in \mathbb{N}$ such that
    $$\{y \in M: \mbox{\em dist} (\sigma_{0,n}(y), \mu) < \epsilon\} \mbox{ \em  is nonempty and open } \  \ \ \forall \ n \geq N.$$
    \end{lemma}

    \noindent{\em Proof: } We take any continuous real function $\varphi \in C^0(M, \mathbb{R})$. By Birkhoff Ergodic Theorem and from the definition of ergodicity (see \cite{libroteoergodica}), we have
    $$\lim_{n \rightarrow + \infty} \frac{1}{n} \sum_{j= 0}^{n-1} \varphi(f^j(y)) = \int \varphi \, d \mu \ \ \mu-\mbox{a.e.} \ y \in M. $$
    From equality (\ref{eqn01}) we obtain
     $$\lim_{n \rightarrow + \infty} \int \varphi \, d\sigma_{0,n}(y)  = \int \varphi \, d \mu \ \ \mu-\mbox{a.e.} \ y \in M. $$
     The last equality holds for all $\varphi \in C^0(M, \mathbb{R})$. So, by the condition (\ref{eqnWeak*Topology}) which defines   the weak$^*$ topology in the space ${\mathcal M}$ of probability measures, we deduce:
     $$\lim_{n \rightarrow + \infty} \sigma_{0,n}  (y) = \mu \ \  \mu-\mbox{a.e.} \ y \in M.$$
     Therefore, for $\mu-$a.e. $y \in M$, for all $\epsilon >0$ there exists $N \geq 1$ such that
     $$\mbox{dist}( \sigma_{0,n} \sigma(y) , \mu) < \epsilon \ \forall \ n \geq N.$$
     We conclude that
     $$\{y \in M: \ \mbox{dist}( \sigma_{0,n} (y) , \mu) < \epsilon \} \neq \emptyset \ \   \forall \ n \geq N. $$

     Now, it is left to prove that, for fixed $\epsilon>0,$ fixed $\mu \in \mathcal M$, and fixed $n \in \mathbb{N}$, the set $\{y \in M: \ \mbox{dist}( \sigma_{0,n} (y) , \mu) < \epsilon \}$ is open in the ambient manifold $M$.
     Since $\{\nu \in\mathcal M: \ \mbox{dist}(\nu, \mu)< \epsilon \}$ is open in the space of probability measures, it is enough to check that the   mapping:
     $$\sigma_{0,n}: M \mapsto {\mathcal M}$$
     is continuous. So, let us prove that  for any convergent sequence $\{x_i\}_{i \in \mathbb{N}} \subset M$, the image sequence $\{\sigma_{0,n}(x_i)\}_{i \in \mathbb{N}}  \subset {\mathcal M}$ converges to $\sigma_{0,n}(x)$ in the weak$^*$ topology, where $$x= \lim x_i \in M.$$

     To apply  condition (\ref{eqnWeak*Topology})
     we consider any continuous real function $\varphi \in C^0(M, \mathbb{R})$. From equality (\ref{eqn01})
     $$\lim_{i \rightarrow + \infty} \int \varphi \, d \sigma_{n,0}(x_i) = \lim_{i \rightarrow + \infty} \frac{1}{n} \sum_{j= 0}^{n-1} \varphi (f^j(x_i)).$$
     Since $\varphi \circ f^j$ is continuous and $\lim x_i= x$ we have
     $$\lim_{i \rightarrow + \infty} \frac{1}{n} \sum_{j= 0}^{n-1} \varphi (f^j(x_i)) = \frac{1}{n} \sum_{j= 0}^{n-1} \varphi (f^j(x)) = \int \varphi \, d \sigma_{n,0}(x).$$
     We deduce that
     $$\lim_{i \rightarrow + \infty} \int \varphi \, d \sigma_{n,0}(x_i)= \int \varphi \, d \sigma_{n,0}(x) \ \ \forall \ \varphi \in C^0(M, \mathbb{R}).$$
     From condition (\ref{eqnWeak*Topology}) we conclude the following equality in the space ${\mathcal M}$ of probability measures:
     $$\lim_{i \rightarrow + \infty} \sigma_{n,0}(x_i) = \sigma_{n,0}(\lim_{i \rightarrow + \infty}x_i), $$
     showing that the mapping $\sigma_{n,0}$ is continuous, and ending the proof of Lemma \ref{lemma0b}.
    \hfill $\Box$

\vspace{.5cm}

  To prove Theorem \ref{theorem1}, we  first   state  the following:
  \begin{lemma}
  \label{lemma1a} Let $f: M \mapsto M$ be continuous, preserve the Lebesgue measure $m$, be ergodic with respect to $m$ and be non-uniquely ergodic.  If $\mu$ is an ergodic probability measure for $f$, then the set  $A_{\mu} $ defined by equality \em  (\ref{eqn11}) \em of  Lemma \em (\ref{lemma0a}) \em has total Lebesgue measure. Thus,  $\mu \in as(x)$ for Lebesgue almost all $x \in M$.
  \end{lemma}

   \noindent{\em Proof: } Denote by $m$ the Lebesgue measure of the manifold $M$, after a rescaling to make $m(M)= 1$.

   From Lemma \ref{lemma0b}, for all $\epsilon >0$ there exists $n_0 \geq 1$ such that
   $$\{y \in M: \ \mbox{dist}(\sigma_{0,n}, \mu) < \epsilon\}   \mbox{ is nonempty and open } \ \ \forall \ n \geq n_0,$$
   Thus $$\bigcup_{n \geq N} \{y \in M: \ \mbox{dist}(\sigma_{0,n}, \mu) < \epsilon\} \mbox{ is nonempty and open }\ \ \forall \ \epsilon>0, \ \ \forall \ N \geq 1. $$
   Define
   $$
   B_N (\epsilon):= \bigcup _{m \geq 0}f^{-m} \Big( \bigcup_{n \geq N}  \{y \in M: \ \mbox{dist}(\sigma_{0,n}, \mu) < \epsilon\} \Big)  = $$ \begin{equation}
   \label{eqn40}\bigcup _{m \geq 0}\bigcup_{n \geq N} f^{-m} \Big (\{y \in M: \ \mbox{dist}(\sigma_{0,n}, \mu) < \epsilon\} \Big).  \end{equation}
   We conclude that $B_N(\epsilon)$ is a nonempty open set in $M$. By construction $f^{-1}(B_N(\epsilon)) \supset B_N(\epsilon)$. Since the Lebesgue measure $m$ is ergodic, we deduce that
   $$m(B_N(\epsilon)) = 0 \mbox{ or } m(B_N(\epsilon)) = 1.$$
   But the set $B_N(\epsilon)$ is nonempty and open, and the Lebesgue measure is positive on nonempty open sets. So
   $$m(B_N(\epsilon))= 1, \mbox{ from where we obtain } m \Big(\bigcap_{N \geq 1} B_N(\epsilon) \Big )=  1  \ \ \forall \ \epsilon >0.$$
So, taking $\epsilon = 1/k, \ \ k \in \mathbb{N}^+$ we deduce that
$$x \in \bigcap_{k \geq 1} \bigcap _{N \geq 1} B_N(1/k)\ \ m-\mbox{a.e. } x \in M.$$
Note that if $0 <\epsilon  < \epsilon'$ then $B_N(\epsilon ) \subset B_N(\epsilon')$, and for any $\epsilon > 0$ there exists $k \in \mathbb{N}^+$ such that $\epsilon < 1/k$. Thus
$$\bigcap_{ \epsilon >0 } \bigcap_{N \geq 1} B_N(\epsilon) \subset \bigcap_{k \geq 1} \bigcap_{N \geq 1} B_N(1/k).$$
But the converse inclusion is obvious because for all $k \in \mathbb{N}^+$, we obtain particular values of $\epsilon = 1/k >0$. Thus
$$\bigcap_{ \epsilon >0 } \bigcap_{N \geq 1} B_N(\epsilon) = \bigcap_{k \geq 1} \bigcap_{N \geq 1} B_N(1/k).$$
 In brief, we have proved that
$$x \in \bigcap_{\epsilon >0} \bigcap _{N \geq 1} B_N(\epsilon)\ \ m-\mbox{a.e. } x \in M.$$
Substituting $B_N$ by its expression in equality (\ref{eqn40}) we conclude:
$$x \in \bigcap_{\epsilon >0} \bigcap _{N \geq 1} \bigcup_{m \geq 0} \bigcup_{n \geq N} f^{-m} \big(\{ y \in M: \ \mbox{dist}(\sigma_{0,n}, \mu) < \epsilon \}  \Big )\ \ m-\mbox{a.e. } x \in M.$$
Finally, applying Lemma \ref{lemma0a} we conclude
$$\mu \in as(x) \ \ m-\mbox{a.e. } x \in M,$$
as wanted. \hfill $\Box$

  \vspace{.5cm}

  \noindent {\bf End of the proof of Theorem \ref{theorem1}}.

  \vspace{.2cm}

  \noindent{\em Proof: } Fix $\epsilon = 1/(2k)$ with $k \in \mathbb{N}^+$. Since the space ${\mathcal M}$ of probability measures is weak$^*$- compact,  the closure $\overline{{\mathcal E}_f}$ is compact. So, there exist a finite covering
  \begin{equation} \label{eqnCovering}{\mathcal C}_k:= \{{\mathcal B}_{1/k}(\mu_{1,k}), \ldots,   {\mathcal B}_{1/k}(\mu_{l(k),k})  \}\end{equation}
  of $\overline{{\mathcal E}_f}$ with open balls ${\mathcal B}_{1/k}(\mu_{i,k}) \subset {\mathcal M}$ of radius $1/k$ and centered in $\mu_{i,k}$. Since the radius $1/k >0$ is fixed, it is not restrictive to take $\mu_{i, k} \in {\mathcal E}_f$ for all $i \in \{1, 2, \ldots, l(k)\}$.

  Since $\mu_{i,k}$ is ergodic for $f$, we can apply Lemma \ref{lemma1a}, to deduce that the following set
  \begin{equation} \label{eqnSetH_k} H_k  := \bigcap _{i= 1 }^{l(k)} \bigcap_{N \geq 1} \bigcup_{m \geq 0} \bigcup_{n \geq N} f^{-m}\Big(\{ y \in M: \ \mbox{dist}(\sigma_{0,n}(y), \mu_{i,k}) < 1/(2k)\}\Big) \subset M \end{equation}
  has total Lebesgue measure.

  Take any $\mu \in {\mathcal E}_f$. Since ${\mathcal C}_k$ covers ${\mathcal E}_f$ there exists $\mu_{i,k}$ such that $\mbox{dist}(\mu, \mu_{i,k}) < 1/(2k)$. Therefore, by the triangle inequality:
$$\{y \in M: \  \mbox{dist}(\sigma_{0,n}(y), \mu_{i,k}) < 1/(2k)\} \subset \{y \in M: \  \mbox{dist}(\sigma_{0,n}(y), \mu  < 1/k\}  \ \forall \ n \in \mathbb{N}.$$

We deduce that
\begin{equation}
\label{eqn50}
H_k \subset  \bigcap_{N \geq 1} \bigcup_{m \geq 0} \bigcup_{n \geq N} f^{-m}\Big(\{ y \in M: \ \mbox{dist}(\sigma_{0,n}(y), \mu)  < 1/k\}\Big)  \ \ \forall \ \mu \in {\mathcal E}_f \ k \in \mathbb{N}^+.\end{equation}
Since $H_k$ has full Lebesgue measure, then the set
$$H = \bigcap_{k \geq 1} H_k$$ also has full Lebesgue measure. From (\ref{eqn50}) we have
 \begin{equation}
 \label{eqnH contained in}
 H \subset \bigcap_{k \geq 1}\bigcap_{N \geq 1}  \bigcup_{m \geq 0} \bigcup_{n \geq N} f^{-m}\Big(\{ y \in M: \ \mbox{dist}(\sigma_{0,n}(y), \mu)  < 1/k\}\Big) \ \ \forall \mu \in {\mathcal E}_f.\end{equation}
We deduce that  for Lebesgue-almost all $x \in M$, the following assertion holds for any ergodic measure $\mu$:
$$ x \in \bigcap_{k \geq 1}\bigcap_{N \geq 1}  \bigcup_{m \geq 0} \bigcup_{n \geq N} f^{-m}\Big(\{ y \in M: \ \mbox{dist}(\sigma_{0,n}(y), \mu)  < 1/k\}\Big) $$
 Since for any $\epsilon >0$ there exists $k \geq 1$ such that $1/k < \epsilon$, we have
$$x \in \bigcap_{\epsilon >0}\bigcap_{N \geq 1}  \bigcup_{m \geq 0} \bigcup_{n \geq N} f^{-m}\Big(\{ y \in M: \ \mbox{dist}(\sigma_{0,n}(y), \mu)  < \epsilon\}\Big).$$
Applying Lemma \ref{lemma0a}, we deduce that $\mu \in as(x)$. We have proved that for Lebesgue-almost all $x \in M$ any ergodic measure $\mu$ belongs to $as(x)$. After Definition \ref{definitionExtemeOscilAE}, the asymptotical statistics of the orbit of $x$ is extremely oscillating. We conclude that    Lebesgue-almost all the orbits exhibit   extremely oscillating asymptotical statistics, as wanted. \hfill $\Box$

  \vspace{.5cm}

  Now, to prove Theorem \ref{theorem0}, we     state the following:

  \begin{lemma}
  \label{lemma1} If $f: M \mapsto M$ is continuous and transitive and if $\mu$ is an ergodic invariant measure for $f$, then the set  $A_{\mu} $ defined in Lemma \em \ref{lemma0a} \em
  is residual. Thus,   $\mu \in as(x)$ for generic $x \in M$.
  \end{lemma}

  \noindent {\em Proof: }

  From Lemma \ref{lemma0b}, for all $\epsilon >0$ the set   $$U_N=\bigcup_{n \geq N} \{y \in M: \ \mbox{dist}(\sigma_{0,n}, \mu) < \epsilon\} \mbox{ is nonempty and open }\ \ \forall \ \epsilon>0, \ \ \forall \ N \geq 1. $$
   Consider the set $B_N(\epsilon) \subset M$ defined by equality
   (\ref{eqn40}):
    $$B_N (\epsilon) = \bigcup_{m \geq 0} f^{-m}(U_N(\epsilon)).$$  Since $f$ is continuous and transitive, from Definition \ref{definitionTransitivo} we obtain that, for any nonempty   open set $V \subset M$, there exists $m \geq 1$ such that $f^m(V) \cap U_N(\epsilon) \neq \emptyset$. In other words, $B_N (\epsilon)$ is dense in $M$. So $B_N$ is open and dense.

So, taking $\epsilon = 1/k, \ \ k \in \mathbb{N}^+$ and applying Definition \ref{definitionGenericOrbit}, we deduce that
$$R:= \bigcap_{k \geq 1} \bigcap _{N \geq 1} B_N(1/k)\ \ \mbox{is residual in }   M.$$
Applying again Definition \ref{definitionGenericOrbit}:
$$x \in \bigcap_{k \geq 1} \bigcap _{N \geq 1} B_N(1/k)\ \ \mbox{ for generic } x \in M.$$
For all $\epsilon >0$ there exists $k \in \mathbb{N}$ such that $\epsilon < 1/k$, and thus
$$B_N (\epsilon) \supset B_N(1/k), \ \ \ \bigcap_{N \geq 1} B_N(\epsilon) \supset R \ \ \forall \ \epsilon >0.$$
We deduce that $$\bigcap_{ \epsilon >0 } \bigcap_{N \geq 1} B_N(\epsilon) \supset R,$$
 is also residual in $M$.
  In other words
$$x \in \bigcap_{\epsilon >0} \bigcap _{N \geq 1} B_N(\epsilon)\ \  \mbox{ for generic } x \in M.$$
Substituting $B_N(\epsilon)$ by its expression in equality (\ref{eqn40}) we conclude:
$$x \in \bigcap_{\epsilon >0} \bigcap _{N \geq 1} \bigcup_{m \geq 0} \bigcup_{n \geq N} f^{-m} \big(\{ y \in M: \ \mbox{dist}(\sigma_{0,n}, \mu) < \epsilon \}  \Big )\ \  \mbox{ for generic } x \in M.$$
Finally, applying Lemma \ref{lemma0a} we conclude
 $\mu \in as(x) \mbox{ for generic } x \in M,$
as wanted. \hfill $\Box$

  \vspace{.5cm}

  \noindent {\bf End of the proof of Theorem \ref{theorem0}}.

  \vspace{.2cm}

  \noindent {\em Proof: }  Fix $\epsilon = 1/(2k)$ with $k \in \mathbb{N}^+$ and construct the finite covering ${\mathcal C}_k$ of $\overline{{\mathcal E}_f}$ by equality (\ref{eqnCovering}), and the set $H_k$ defined by equality (\ref{eqnSetH_k}).
  Since the measures $\mu_{i,k}$ are ergodic for $f$, we can apply Lemma \ref{lemma1}, to deduce that the set $H_k$ is residual for all $k \in \mathbb{N}^+$.

  Take any $\mu \in {\mathcal E}_f$. Since ${\mathcal C}_k$ covers ${\mathcal E}_f$ there exists $\mu_{i,k}$ such that $\mbox{dist}(\mu, \mu_{i,k}) < 1/(2k)$. Therefore, by the triangle inequality, we deduce assertion (\ref{eqn50}).

Since $H_k$ is residual, then the set
 $H = \bigcap_{k \geq 1} H_k$  is also residual. From (\ref{eqn50}) we have $$ H \subset \bigcap_{k \geq 1}\bigcap_{N \geq 1}  \bigcup_{m \geq 0} \bigcup_{n \geq N} f^{-m}\Big(\{ y \in M: \ \mbox{dist}(\sigma_{0,n}(y), \mu)  < 1/k\}\Big) \ \ \forall \mu \in {\mathcal E}_f.$$
We deduce that  for generic $x \in M$, the following assertion holds for any ergodic measure $\mu$:
$$ x \in \bigcap_{k \geq 1}\bigcap_{N \geq 1}  \bigcup_{m \geq 0} \bigcup_{n \geq N} f^{-m}\Big(\{ y \in M: \ \mbox{dist}(\sigma_{0,n}(y), \mu)  < 1/k\}\Big) $$
 Since for any $\epsilon >0$ there exists $k \geq 1$ such that $1/k < \epsilon$, we have
$$x \in \bigcap_{\epsilon >0}\bigcap_{N \geq 1}  \bigcup_{m \geq 0} \bigcup_{n \geq N} f^{-m}\Big(\{ y \in M: \ \mbox{dist}(\sigma_{0,n}(y), \mu)  < \epsilon\}\Big).$$
Applying Lemma \ref{lemma0a}, we deduce that $\mu \in as(x)$. We have proved that for generic $x \in M$ any ergodic measure $\mu$ belongs to $as(x)$. After Definition \ref{definitionExtemeOscilAE}, the asymptotical statistics of the orbit of $x$ is extremely oscillating. We conclude that    the generic orbits of $f$ exhibit   extremely oscillating asymptotical statistics, as wanted. \hfill $\Box$

\vspace{.8cm}

\noindent{\bf Acknowledgements: } \thanks{The author thanks the Editor and the anonymous Referee. She thanks the partial support of
 \lq\lq Agencia Nacional de Investigaci\'{o}n e Innovaci\'{o}n\rq\rq \ (ANII), \lq\lq Comisi\'{o}n Sectorial de Investigaci\'{o}n Cient\'{\i}fica\rq\rq \ (CSIC) of \lq\lq Universidad de la Rep\'{u}blica\rq\rq, and \lq\lq Premio L'Or\'{e}al-UNESCO-DICYT \rq\rq (the three institutions of Uruguay).
}


\end{document}